\documentclass[letter,11pt]{article}
\usepackage{amsmath}
\usepackage{bbding}
\usepackage{amssymb}
\usepackage[margin=1in]{geometry}
\usepackage{graphics,graphicx}
\usepackage[round]{natbib}

\usepackage{amsmath}
\usepackage{bbding}
\usepackage{amssymb}
\usepackage{lineno,hyperref}
\usepackage{amssymb,amsmath,amsthm}
\usepackage{color}
\usepackage{graphicx}
\usepackage{subfigure}
\usepackage{epstopdf}
\usepackage{graphics,epsfig}
\usepackage{algorithm}

\newcommand{\argmin}{\arg\!\min}

\usepackage{color}

\begin{document}
\bibliographystyle{plainnat} 

\title{Estimation of CD4+ T Cell Count Parameters in HIV/AIDS Patients Based on Real-time Nonlinear Receding Horizon Control}

\author{Fei Sun \thanks{Post-Doctoral Researcher, Control Science and Dynamic Systems (CSDy) Laboratory, Aerospace Engineering Department, San Jose State University, fei.sun@sjsu.edu} and Kamran Turkoglu \thanks{Assistant Professor, Aerospace Engineering Department, San Jose State University, kamran.turkoglu@sjsu.edu} \\
Aerospace Engineering, San Jose State University, \\ San Jose 95192 CA}
 
\maketitle
 
\begin{abstract}
An increasing number of control techniques are introduced to HIV infection problem to explore the options of helping clinical testing, optimizing drug treatments and to study the drug resistance situations. In such cases, complete/accurate knowledge of the HIV model and/or parameters is critical not only to monitor the dynamics of the system, but also to adjust the therapy accordingly. In those studies, existence of any type of unknown parameters imposes severe set-backs and becomes problematic for the treatment of the patients. In this work, we develop a real-time adaptive nonlinear receding horizon control approach to aid such scenarios and to estimate unknown constant/time-varying parameters of nonlinear HIV system models. For this purpose, the problem of estimation is updated by a series of finite-time optimization problem which can be solved by backwards sweep Riccati method in real time without employing any iteration techniques. The simulation results demonstrates the fact that proposed algorithm is able to estimate unknown constant/time-varying parameters of HIV/AIDS model effectively and provide a unique, adaptive solution methodology to an important open problem.
\end{abstract}

\section{Introduction: Background and Motivation}
When first discovered by Gallo and his team (\cite{gallo1984,popovic1984,sarngadharan1984,schupbach1984}), the working principles of Human Immunodeficiency Virus (HIV) were not well known. Throughout the years, many scientists from different fields, ranging from medicine, bio-chemistry to mechanics, reported valuable insights regarding the dynamics of HIV. While, some of these models are purely theory inspired, some are clinical trial (and/or statistics) based. But one common goal in majority of those studies is to identify the working principles of Human Immunodeficiency Virus (HIV). As will be discussed in the following sections in more detail, commonly taken modeling approach in HIV architecture is to define healthy and infected CD4+ T cells together with the viron dynamics, as reported in \cite{xia2003estimation}. In some of the studies in existing literature, researchers also introduce extra components in HIV and cell dynamics (such as mutant cells) to reflect the true nature of HIV dynamics, as stated in \cite{kim2015}.

In existing work on HIV modeling and control, there are many valuable studies concentrating on the use of modern and advanced control theory, as a brief survey by \cite{perelson1999mathematical} outlines some of the earlier efforts. In modeling and control part, there is a common branching between the choice of models and estimation routines used. While some earlier models (like \cite{perelson1996hiv,wu1999population}) discussed mostly statistics based approach, there are some models that incorporate linear dynamics (\cite{lemos2012control}). It is important to note that linear models provide very limited and constrained information, thus in recent years there is more effort in the scientific community to incorporate the true non-linear and highly coupled nature of the HIV infection problem (\cite{david2011receding}).

In one of those studies \cite{lemos2012control} proposes feedback linearization and LQ regulation based linear state feedback strategies. In another recent study, \cite{kim2015adaptive} demonstrate a $\sigma-modification$ based adaptive estimation routine of HIV parameters in presence of mutant cells. \cite{david2011receding} discuss a methodology based on an iterative receding horizon control, kalman filtering and state estimation framework for HIV. \cite{liang2010estimation} presents a multi-stage smoothing-based (MSSB) and the spline-enhanced nonlinear least squares (SNLS) approach, to estimate all HIV viral dynamic parameters. \cite{pinheiro2011}, investigates a nonlinear model predictive control algorithm embedding nonlinear multirate state estimation with an extended Kalman filter for minimized drug delivery rate for effective treatment. \cite{casal2013multiple} implements an estimation algorithm using a multiple-model adaptive estimation approach with a bank of moving horizon estimators on HIV parameter estimation dynamics. \cite{hartmann2012online}, in their study, explore nonlinear bayesian filtering through extended kalman filter and particle filter for joint state and time-varying parameter estimation HIV problem. \cite{david2011receding} discuss a receding horizon based methodology for the estimation problem and explore how ideas from stochastic estimation can be used to create a better treatment methodology.

As an extension to existing studies on estimation, there are also reported efforts (such as \cite{portelo2013parameter,wu2008parameter}) where the emphasis was given on identifiability aspect of the parameters and parameter estimation routine of HIV problem(s).

After so many years, it still remains as a challenging problem to establish a routine for unknown parameter estimation in HIV related studies, both from clinical and/or mathematical stand point. While the emphasis of most HIV studies were/are two fold (namely modeling and estimation), the main concentration of this study remains to be on parameter estimation and control system design for HIV infected cell models. We will use presented models to not only estimate unknown parameters, but also will we design non-linear control schemes to aid the control of infected cells. This architecture is different than the existing studies because of the following factors:

\begin{itemize}
 \item In this study we provide estimates to some of HIV infection model constant/time-varying parameters based on the adaptive receding horizon control methodology without the knowledge of the true value of parameters. With the adaptive identifier, the estimation problem is converted into a series of finite-time optimization problem which can be presented as nonlinear Riccati equations solving directly in real time by the backward sweep method without employing any iteration techniques.
 \item In addition, in presence of certain assumptions, we also demonstrate the fact that when such  parameter estimation routine exists, obtained results present better control of HIV structure in simulations, in combination with nonlinear receding horizon control architecture, which is the second contribution to the existing literature.
\end{itemize}

Here, in this study, it is assumed that the states are available to be used in measurements and we discuss further details under the following organization scheme of the paper: In Section-\ref{sec:prob_formulation}, we analyze the model used in estimation of HIV parameters. In Section- \ref{sec:adaptiv_nrhc}, we present adaptive nonlinear receding horizon control based methodology, where in Section-\ref{sec:backward_sweep} backward sweep dynamics are derived. Section-\ref{sec:convergence} investigates convergence analysis of proposed methodology. In Section-\ref{sec:sims} results and outcomes (with two sample scenarios of the study) are iterated in further detail, and with Section-\ref{sec:conclusions} the paper is concluded. 

%
%

\section{The HIV Model}\label{sec:prob_formulation}
\vspace{-2pt}
In this paper, we consider the following set of differential equations as the HIV infection model, which is inherited from \cite{xia2003estimation} and \cite{liang2010}:

\begin{equation}\label{eq:hiv_system}
\begin{cases}
\dot{x}_1=s-dx_1-\beta x_1x_3,\\
\dot{x}_2=\beta x_1x_3-{\mu}_1x_2,\\
\dot{x}_3=kx_2-{\mu}_2x_3,
\end{cases}
\end{equation}
Here, state variables are $x_1$, $x_2$ and $x_3$, and dynamic parameters are $s,d,\beta,\mu_1,k,\mu_2$. In this set-up, $x_1$ represents the concentration of uninfected CD4+T cells in patient blood per cubic millimeter in presence of HIV, $x_2$ is the concentration of infected cells and $x_3$ is the concentration of free virus particles. $s$ denotes the proliferation rate of uninfected cells, $d$ denotes the death rate of uninfected cells, $\beta$ represents the infection rate, $\mu_1$ is the death rate of infected cells, $k$ is the production rate of free virons, and $\mu_2$ is the death rate of free virus particles. 

In Eq.(\ref{eq:hiv_system}), first equation describes the population dynamics of uninfected CD4+T cells which are created at rate $s$ from sources within the body, died at rate $d$ and infected by HIV at rate $\beta$. Here, the infection is represented by a "mass-action" term, namely $\beta x_1x_3$. The second equation defines the population dynamics of infected cells which are generated at rate $\beta$ and died at rate $\mu_1$. The last equation defines the dynamics of concentration of free virus particles which are produced by infected cells at rate $k$ and died at rate $\mu_2$.

With this formulation, the goal of this work is to estimate dynamic parameters for more effective drug treatment. Based on previously conducted clinical studies (\cite{liang2010}), here we consider two options for parameters of interest: (i) they could be constant or, (ii) they could be time-varying. In practice, some parameters can be fixed and only the remaining parameters are needed to be estimated, while all those parameters could be constant or time-varying. From the clinical standpoint, all the above mentioned three variables which are the concentration of healthy cells-$x_1$, the concentration of infected cells-$x_2$ and the load of free virus-$x_3$ can be measured.

\section{Adaptive Nonlinear Receding Horizon Control Methodology}\label{sec:adaptiv_nrhc}

In order to estimate desired unknown parameters, the HIV model (presented in Eq.(\ref{eq:hiv_system})) is modified as the following nonlinear system of the form:
\begin{equation}\label{eq:sys_hiv}
\dot{x}=f(x(t),\Theta(t))=g(x)+D(x)\Theta(t),\\
\end{equation}
where $x\in R^{n}$ is the state vector, $A,B\in R^{n\times n}$ is the the linear coefficient
matrix and $g(x):R^{n}\rightarrow R^{n}$ is the nonlinear part of system in Eq.\eqref{eq:sys_hiv}. $D(x):R^{n}\rightarrow R^{n\times p}$ is a known function vector and $\Theta(t) \in R^{p}$ denotes the unknown dynamic parameters. 

In order to identify the unknown \emph{dynamic} parameters, we are interested in constructing a drive-response scheme. We take the HIV model in Eq. \eqref{eq:sys_hiv} as the drive system. And the corresponding controlled response system is designed as follows

\begin{equation}\label{eq:sys_hiv_response}
\dot{y}=f(y(t),\hat{\Theta}(t),u(t))=g(y)+D(y)\hat{\Theta}(t)+u(t),\\
\end{equation}
where $y\in R^{n}$ is the state vector, $\hat{\Theta}(t)$ represents the estimated parameters which can be constant or time-varying and $u(t)$ is the controller. Here, function $g(\cdot)$ and $D(\cdot)$ satisfy the global Lipschitz condition. Therefore there exist positive constants $\beta_1$ and $\beta_2$ such that
\begin{equation*}
\begin{split}
\|g(y)-g(x)\|\le&\beta_1\|y-x\|,\\
\|D(y)-D(x)\|\le&\beta_2\|y-x\|.
\end{split}
\end{equation*}
This condition is satisfied if the Jacobians $\frac{\partial g}{\partial x}$, $\frac{\partial g}{\partial y}$, $\frac{\partial D}{\partial x}$ and $\frac{\partial D}{\partial y}$ are uniformly bounded.

The synchronization error between the drive and the response system becomes an important part of this analysis, and in this study it is defined as
$$e(t) = y(t) -x(t).$$
Here, the parameter updating law is designed by the following standard adaptive estimator
\begin{equation}\label{estimator}
\dot{\hat{\Theta}}(t)=-D(y)^Te(t),
\end{equation}
where the superscript $T$ is denoted as the matrix transposition. 

\textbf{Assumption-1:} $D(x)$ satisfies the persistent excitation condition, which means there exists strictly positive constants $\alpha_1$ and $\delta_1$ such that for any $t\ge 0$,
$$\int_t^{t+\delta_1}D(y(\tau))D(y(\tau))^Td\tau\ge\alpha_1I,$$
where $I$ is the identity matrix with the same dimension of $D(x)$.

In such formulation,  we combine the adaptive estimator with the real-time nonlinear receding horizon control to identify the unknown parameters as well as generating the control protocols. The following optimization problem is utilized to build the controller:

\textbf{Problem-1:}
\begin{equation}
u^*(t)=\argmin_{u(t)}J(y(t),x(t),u(t),\hat{\theta}(t))
\end{equation}
\hskip 70pt subject to
\begin{equation*}
\dot{y}=g(y)+D(y)\hat{\Theta}(t)+u(t),
\end{equation*}
where the performance index associated with the synchronization error and the controller is designed as follows:

\begin{equation}\label{eq:sync_cost_func}
\begin{split}
J=&\int_t^{t+T}L[e(\tau),u(\tau)]{\rm d}\tau,\\
=&\int_t^{t+T}[e(t)^TQe(t)+u(t)^TRu(t)]{\rm d}\tau,
\end{split}
\end{equation}
where weighting matrices $Q$ and $R$ are positive definite. With the construction of the synchronization cost function given in Eq.\eqref{eq:sync_cost_func}, the estimated parameters are updated and the optimal control strategy are generated simultaneously. For this purpose, we utilize the powerful nature of \emph{real-time nonlinear receding horizon control algorithm} to minimize the associated synchronization cost function. The corresponding performance index evaluates the performance from the present time $t$ to the finite future $t+T$, and is minimized for each time segment $t$ starting from $e(t)$. With this structure, it is possible to convert the present receding horizon control problem into a family of finite horizon optimal control problems on the artificial $\tau$ axis parameterized by time $t$.

According to the well-known first order necessary conditions of optimality, we can obtain the two-point boundary value problem (TPBVP) \cite{bryson1975}, as follows:

\begin{equation}\label{eq:tpbvp}
\begin{split}
y_\tau^* = H_{\lambda}^T , \hskip 20pt y^* = &y(t), \hskip 20pt \lambda^*(\tau,t) = -H_y ^T \\
 \lambda^*(T,t) =&~ 0, \hskip 20pt  H_{u} = 0
\end{split}
\end{equation}
where $\lambda(\tau,t)\in R^{n}$ denotes the co-state and $y_\tau$ denotes the partial derivative of $y$ with respect to $\tau$, and so on. $H$ is the Hamiltonian defined of the form
\begin{equation}\label{eq:tpbvp_ham}
\begin{split}
H &= L + \lambda^{*T}\dot{y}\\
&=(e^TQe+u^TRu)+ \lambda^{*T}[Ay+g(y)+D(y)\hat{\Theta}(t)].
\end{split}
\end{equation}
In Eqs.\eqref{eq:tpbvp}-\eqref{eq:tpbvp_ham}, $(~~)^*$ denotes a variable in the optimal control problem so as to distinguish it from its original term.

Using this formulation, the optimal controller can be calculated as
\begin{equation}\label{ctrl}
u(t) = arg\{H_u[y(t),\lambda(t),\hat{\Theta}(t),x(t)] = 0\}
\end{equation}

In order to solve the above equation of $u(t)$ in real-time, the TPBVP is to be rewritten as another form. Since the state and co-state at $\tau=T$ are determined by the TPBVP in Eq.\eqref{eq:tpbvp} from the state and co-state at $\tau=0$, the TPBVP can be considered as a nonlinear algebraic equation regarding to the co-state at $\tau=0$ as
\begin{equation}
F(\lambda(t),y(t),T,t)=\lambda^*(T,t)=0.
\end{equation}

Due to that the nonlinear equation $F(\lambda(t),y(t),T,t)$ has to be satisfied at any time $t$, $\frac{{\rm d}F}{{\rm d}t}=0$ holds along the trajectory of the closed-loop system of the receding horizon control. If $T$ is a smooth function of time $t$, the solution of $F(\lambda(t),y(t),T(t),t)$ can be tracked with respect to time $t$. The optimal controller can be computed from  Eq.\eqref{ctrl} based on $y(t), \lambda(t), \hat{\Theta}(t), x(t)$ where the ordinary differential equation of $\lambda(t)$ can be solved numerically, in real-time, without any need of an iterative optimization routine. Meanwhile, the estimated  parameters $\hat{\Theta}(t)$ are updated from Eq.\eqref{estimator} based on current $x(t), y(t), u(t)$ accordingly. However, in practice, numerical errors associated with the solution may accumulate as the integration proceeds and therefore we employ some correction techniques to correct such errors in the solution. A stabilized continuation method \cite{kabamba1987,ohtsuka1994,ohtsuka1997,ohtsuka1998} is applied here. According to this method, it is possible to rewrite the statement as
\begin{equation}\label{cm}
\frac{{\rm d}F}{{\rm d}t}=-A_sF,
\end{equation}
where $A_s$ denotes a stable matrix to make the solution converge to zero asymptotically, so that
\begin{equation}
F=F_0e^{-A_st}\rightarrow 0 \,(t\rightarrow \infty).
\end{equation}
Here, $F_0$ denotes the initial value of $F$, which shows that the numerical errors will eventually be attenuated through the integration process.

\section{Backward-sweep Method:}\label{sec:backward_sweep}

In order to compute the controller $u(t)$ from Eq.\eqref{ctrl}, we first
integrate the differential equation of $\lambda(t)$ in real time. The partial differentation of Eqs.\eqref{eq:tpbvp} with respect to time $t$ and $\tau$ can be written as,
\begin{equation}
\begin{split}
&\delta \dot{y}=f_y\delta y+f_{u}\delta u,\\
&\delta \dot{\lambda}=-H_{yy}\delta y-H_{y\lambda}\delta \lambda-H_{yu}\delta u\\
&0=H_{uy}\delta y+f_{u}^T\delta \lambda+H_{uu}\delta u.
\end{split}
\end{equation}

Since $\delta u=-H^{-1}_{uu}(H_{uy}\delta y+f_{u}^T\delta \lambda$), we have
\begin{align*}
&\delta \dot{y}=(f_y-f_{u}H^{-1}_{uu}H_{uy})\delta y-f_{u}H^{-1}_{uu}f_{u}^T\delta \lambda,\\
&\delta \dot{\lambda}=-(H_{yy}-H_{yu}H^{-1}_{uu}H_{uy})\delta y-(f_{y}^T-H_{yu}H^{-1}_{uu}f_{u}^T)\delta \lambda,
\end{align*} 
which converts the problem to the following linear differential equation:
\begin{equation}\label{pd}
\frac {\partial}{\partial \tau}\begin{bmatrix}y^*_t-y^*_{\tau}\\
\lambda^*_t-\lambda^*_{\tau} \end{bmatrix}=\begin{bmatrix}G&-L\\
-K&-G^T \end{bmatrix}\begin{bmatrix}y^*_t-y^*_{\tau}\\
\lambda^*_t-\lambda^*_{\tau} \end{bmatrix}
\end{equation}
where $G=f_y-f_{u}H^{-1}_{uu}H_{u y}$, $L=f_{u}H^{-1}_{uu}f_{u}^T$, $K=H_{yy}-H_{yu}H^{-1}_{uu}H_{u y}$. And the matrix $H_{uu}$ should be nonsingular.

The derivative of the nonlinear function $F$ with respect to time $t$ is given by
\begin{equation}\label{df}
\frac{{\rm d}F}{{\rm d}t}=\lambda^*_t(T,t)+\lambda^*_{\tau}(T,t)\frac{{\rm d}T}{{\rm d}t}.
\end{equation}

In order to reduce the computational cost without resorting to any approximation technique, the backward-sweep method is implemented here. The relationship between the co-state and other variables is then expressed as the followings:
\begin{equation}\label{relation}
\lambda^*_t-\lambda^*_\tau=S(\tau,t)(y^*_t-y^*_\tau)+c(\tau,t),
\end{equation}
where
\begin{equation}\label{sc}
\begin{split}
S_{\tau}&=-G^TS-SG+SLS-K, \\
c_{\tau}&=-(G^T-SL)c,
\end{split}
\end{equation}
and
\begin{equation}\label{sct}
\begin{split}
S(T,t)&=0,\\
c(T,t)&=H_y^T\mid_{\tau=T}(1+\frac{{\rm d}T}{{\rm d}t})-A_sF.
\end{split}
\end{equation}
With this formulation, the differential equation of $\lambda(t)$ to be integrated in real time becomes:
\begin{align}\label{dl}
\frac{{\rm d}\lambda(t)}{{\rm d}t}=-H_y^T+c(0,t).
\end{align}

Here, at each time $t$, the Euler-Lagrange equations (Eqs.\eqref{eq:tpbvp}) are integrated forward along the artificial $\tau$ axis (which depends on the designated horizon time), and then Eqs.\eqref{sc} are integrated backward with terminal conditions provided in Eq.\eqref{sct}. Next, the differential equation of $\lambda(t)$ is integrated for one step along the $t$ axis so as to compute the optimal control policies from Eq.\eqref{ctrl}. If the matrix $H_{uu}$ is nonsingular, the algorithm is executable regardless of controllability or stabilizability or the system.

The adaptive NRHC method for computing the optimal control $u(t)$ as well as estimating the unknown parameters are summarized in Algorithm-1, where $t_s >0$ denotes the sampling time.

\begin{algorithm}\label{alg_1}
 \caption{ $u^*(t)=\argmin_{u(t)}J(y(t),x(t),u(t),\hat{\theta}(t))$} 
 (1) Set $t=0$ and initial state $y=y(0)$.\\
 (2) For $t^{\prime}\in[t,t+t_s]$, integrate the defined TPBVP in \eqref{eq:tpbvp} forward from $t$ to $t+T$, then integrate \eqref{sc} backward with terminal conditions provided in \eqref{sct} from $t+T$ to $t$. \\
 (3) Integrate the differential equation of $\lambda(t)$, from $t$ to $t+ t_s$.\\
 (4) At time $t+t_s$, compute $u^*$ by Eq.\eqref{ctrl} with the terminal values of $y(t)$, $\lambda(t)$, $\hat{\theta}(t)$, $x(t)$ and update the estimated parameters $\hat{\theta}(t)$ by Eq.\eqref{estimator}.\\
 (5) Set $t=t+t_s$, return to Step-(2).
\end{algorithm}

\section{Convergence and Stability Analysis of NRHC Estimation Problem}\label{sec:convergence}

In order to ensure the closed-loop stability of the proposed adaptive nonlinear receding horizon control scheme, we first introduce some definitions.

In this regard, we assume the sublevel sets 
\begin{equation*}
\Gamma_r^{\infty}=\{e\in\Gamma^{\infty}:J^*_{\infty}<r^2\}
\end{equation*}
are compact and path connected where $J^*_{\infty}=\int_0^{\infty}[e^{*T}Qe^*+u^{*T}Ru^*]{\rm d}\tau$ and moreover $\Gamma^{\infty}=\cup_{r\ge 0}\Gamma_r^{\infty}$. We use $r^2$ here to reflect the fact that the cost function is quadratically bounded. And therefore the sublevel set of $\Gamma_r^{T}=\{e\in\Gamma^{\infty}:J^*_{T}<r^2\}$ where $J^*_{T}=\int_t^{t+T}[e^{*T}Qe^*+u^{*T}Ru^*]{\rm d}\tau$.

\textbf{Lemma-1:} (Dini \cite{jadbabaie2005} ) Let $\{f_n\}$ be a sequence of upper semi-continuous, real-valued functions on a countably compact space $X$, and suppose that for each $x\in X$, the sequence $\{f_n(x)\}$ decreases monotonically to zero. Then the convergence is uniform.

\textbf{Theorem-1:} \cite{jadbabaie2005} Let $r$ be given as $r>0$ and suppose that the terminal cost is equal to zero. For each sampling time $\Delta >0$, there exists a horizon window  $T^*<\infty$ such that, for any $T \geqslant T^*$, the receding horizon scheme is asymptotically stabilizing. 

\begin{proof}
By the principle of optimality, we have
\begin{equation*}
J^*_T(e)=\int_{t}^{t+\Delta}(e_T^{*T}Qe_T^*+u_T^{*T} Ru_T^*)d\tau+J^*_{T-\Delta}(e_T^*)
\end{equation*}
where $\Delta \in [t,t+T]$ is the sampling time and $J^*_{\Delta}(e)=\int_{t}^{t+\Delta}(e_T^{*T}Qe_T^*+u_T^{*T} Ru_T^*)d\tau$, so that
\begin{equation*}
\begin{split}
J^*_{T-\Delta}(e_T^*)-J^*_{T-\Delta}(e)&=J^*_{T}(e)-J^*_{T-\Delta}(e)-\int_t^{t+\Delta}(e_T^{*T}Qe_T^{*T}+u_T^{*T} Ru_T^{*T})d\tau\\
&\leq J^*_{\Delta}(e)+J^*_{T}(e)-J^*_{T-\Delta}(e)
\end{split}
\end{equation*}
Since the terminal cost is equal to zero, it is clear that $T_1<T_2$. This implies that $J^*_{T_1}(e)<J^*_{T_2}(e)$ holds for all $e$ so that
\begin{equation*}
J^*_{T-\Delta}(e_T^*)-J^*_{T-\Delta}(e) \leq J^*_{\Delta}(e)+J^*_{\infty}(e)-J^*_{T-\Delta}(e).
\end{equation*}
is satisfied.
If we can show, for example, that there exists a $T^*$ such that $T>T^*$ yields into
\begin{equation*}
J^*_{\infty}(e)-J^*_{T-\Delta}(e) \leq \frac{1}{2}J^*_{\Delta}(e)
\end{equation*}
for all $e \in \Gamma_r^{\infty}$, stability over any sublevel set of $J^*_{T-\Delta}(\cdot)$ that is contained in $\Gamma_r^{\infty}$ will be assured. To that end, define, for $e \in \Gamma_r^{\infty}$

\begin{equation*}
\psi_T(e)=
\begin{cases}
\frac{J^*_{\infty}(e)-J^*_{T-\Delta}(e)}{J^*_{\delta}(e)}, & \quad e \neq 0 \\
\lim \sup_{x\rightarrow 0} \psi_T(e), & \quad e=0
\end{cases}
\end{equation*}
where $\psi_T(\cdot)$ is upper semicontinuous on $\Gamma_r^{\infty}$. It is clear that $\psi_T(\cdot)$ is a monotonically decreasing family of upper semicontinuous functions defined over the compact set $\Gamma_r^{\infty}$. Thus, by Dini's theorem (as stated in \cite{jadbabaie2005}), there exists a $T^*<\infty$ such that $\psi_T(e)<\frac{1}{2}$ for all $e \in \Gamma_r^{\infty}$ and all $T\ge T^*$. Here, for each $r_1>0$ we have $\Gamma_{r_1}^{T-\Delta}\subset \Gamma_r^{\infty}$ satisfied, leading to

\begin{equation*}
J^*_{T-\Delta}(e_T^*)-J^*_{T-\Delta}(e) \leq -\frac{1}{2}J^*_{\Delta}(e)
\end{equation*}
for all $e \in \Gamma_{r_1}^{T-\Delta}$.

\end{proof}

\textbf{Corollary-1:} Consider the nonlinear HIV system given in Eq.\eqref{eq:hiv_system} and the system satisfies the \textbf{Assumption-1}. For the proposed control protocol in \textbf{Problem-1} and adaptive parameter estimator, based on \textbf{Theorem-1}, there exists a large enough value of horizon $T$ which guarantees the synchronous error $e$ to remain asymptotically stable to achieve synchronization and estimate unknown parameters in the HIV system.

With this result, when the optimization horizon is chosen to be sufficiently long, the non-increasing monotonicity of the cost function becomes a sufficient condition for the stability and the unknown parameters can be estimated through this process, without any given reference trajectory for the parameters.

\section{Simulation Studies}\label{sec:sims}

In this section, we provide concrete examples to validate the performance of the proposed adaptive nonlinear receding horizon control method on the HIV model in Eq. \eqref{eq:hiv_system} through two simulations. 
\subsection{Case-1: Constant parameters}
We first consider the case that all parameters are constant. The following parameter values are given to generate the simulation data:$s=36$, $d=0.108$, $\beta=9\times 10^{-5}$, ${\mu}_1=0.5$, $k=500$, ${\mu}_2=3$. We assume that $s$, $\mu_1$ and $k$ are unknown parameters. The initial conditions of state variables are set to be $x_1(0)=1000$, $x_2(0)=10$, $x_3(0)=1000$, and $y_1(0)=200$, $y_2(0)=50$, $y_3(0)=20000$. And the initial conditions of estimated parameters are set to be $s(0)=1$, $\mu_1=1$ and $k=1$.

The weighting matrices in the performance index are designed as $Q=R=\text{diag}(1,1,1)$. And the horizon $T$ in the performance index is given by
\begin{equation}\label{ex_h}
T(t)=T_f(1-e^{-\alpha t}),
\end{equation}
where $T_f=0.1$ and $\alpha=0.01$. And the stable matrix in Eq. \eqref{cm} $A_s=60I$.

The simulation is implemented in MATLAB, where the sampling time $t_s$ is $0.01$ day and the time step on the artificial $\tau$ axis is $0.005$ day. Fig. \ref{fig_sys} depicts the trajectories of the drive HIV system and the response HIV system which shows the response system tracks the drive system after $7$ days under the proposed adaptive nonlinear receding horizon control methodology. Fig. \ref{fig_u} describes the optimal control strategy generated by Algorithm-1 and Fig. \ref{fig_par} describe the true values of system parameters and the estimated results which converge to the true values after $7$ days. These figures show the optimal control policies result in the system synchronization and the parameter convergence simultaneously, which clearly demonstrates the effectiveness of the proposed algorithm.

 \begin{figure}
 \centering
 \vspace*{-1.1cm}
    \includegraphics[width=13.5cm]{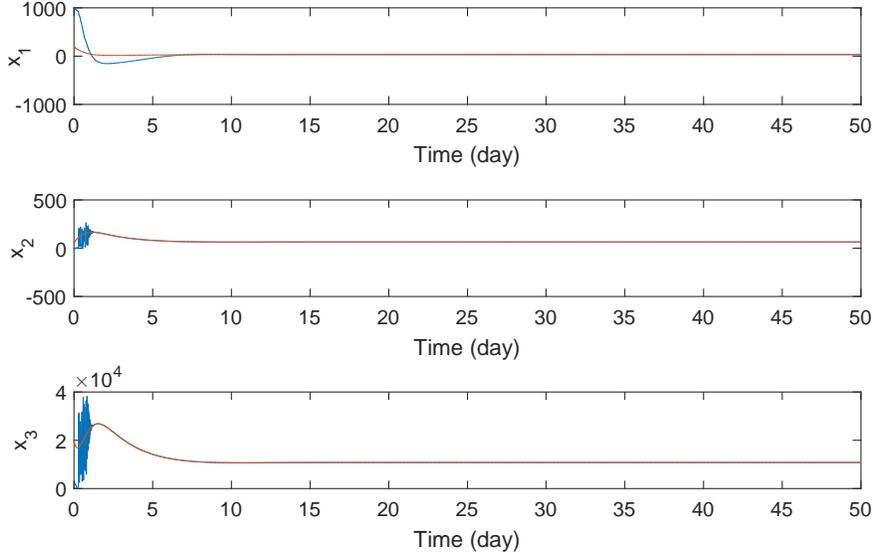}
    \vspace*{-0.5cm}
    \caption{\label{fig_sys}  (color online) Trajectories of three states of HIV model in Eq. \eqref{hiv_system} where the red lines represent the dynamical evolution of drive system and the blue lines represent the evolution of response system with unknown parameters.}
 \end{figure}

\begin{figure}
 \centering
 \vspace*{-1.1cm}
    \includegraphics[width=13.5cm]{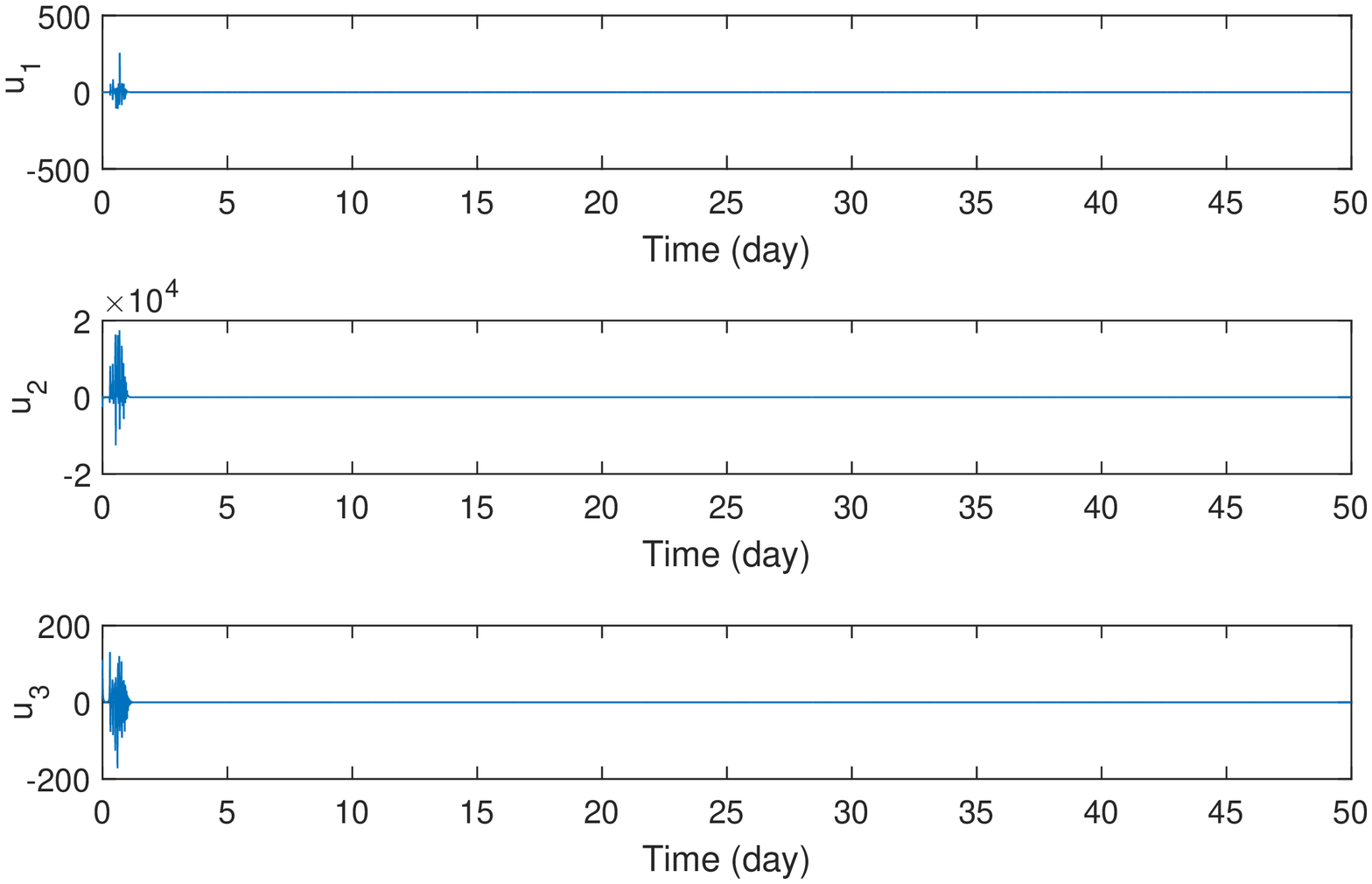}
    \vspace*{-0.5cm}
    \caption{\label{fig_u}  (color online) Trajectories of optimal control strategies generated by adaptive nonlinear receding horizon control method in Algorithm-1.}
 \end{figure}
 
 \begin{figure}
 \centering
 \vspace*{-1.1cm}
    \includegraphics[width=13.5cm]{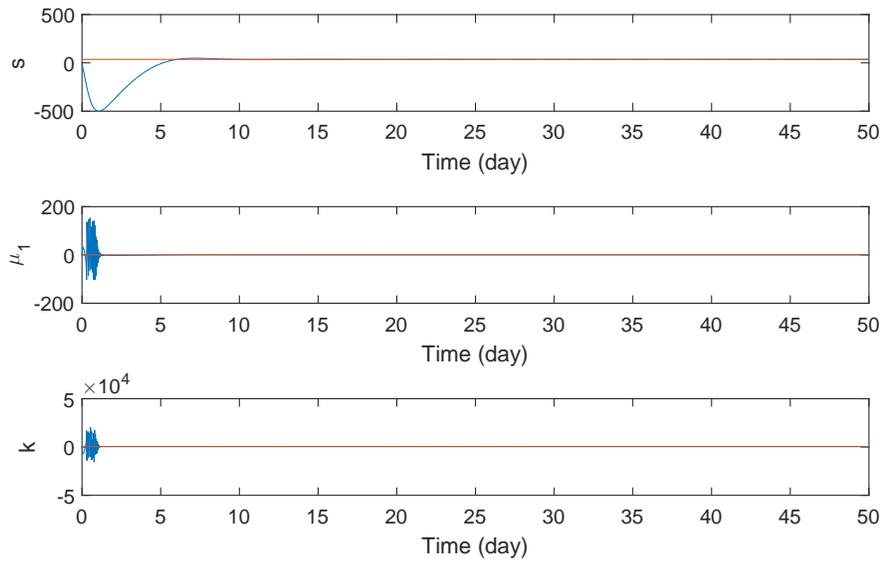}
    \vspace*{-0.5cm}
    \caption{\label{fig_par}  (color online) Trajectories of estimated dynamic parameters of HIV model where the red lines denote their true values and the blue lines denote the estimated results.}
 \end{figure}

 \subsection{Case-2: Time-varying parameters}
 
Next we consider the case where there is a time-varying parameter in the system. The following parameter values are given to generate the simulation data:$s=36\times(1-0.9\cos(\pi t/1000))$, $d=0.108$, $\beta=9\times 10^{-5}$, ${\mu}_1=0.5$, $k=500$, ${\mu}_2=3$. $s$, $\mu_1$ and $k$ are unknown parameters to be estimated. The initial conditions of state variables are the same as previous simulation where $x_1(0)=1000$, $x_2(0)=10$, $x_3(0)=1000$, and $y_1(0)=200$, $y_2(0)=50$, $y_3(0)=20000$. And the initial conditions of estimated parameters are also set to be $s(0)=1$, $\mu_1=1$ and $k=1$.

The weighting matrices in the performance index are also designed as $Q=R=\text{diag}(1,1,1)$. And the horizon $T$ in the performance index is given in Eq. \eqref{ex_h} where $T_f=0.1$ and $\alpha=0.01$. And the stable matrix $A_s=60I$.

The simulation is implemented in MATLAB, where the sampling time $t_s$ is $0.01$ day and the time step on the artificial $\tau$ axis is $0.005$ day. Fig. \ref{fig_sys} depicts the trajectories of the drive HIV system and the response HIV system with time-varying parameter which shows the response system tracks the drive system after around $50$ days under the proposed adaptive nonlinear receding horizon control methodology. Fig. \ref{fig_u} describes the optimal control strategy generated by Algorithm-1 and Fig. \ref{fig_par} describe the true values of system constant and time-varying parameters and the corresponding estimated results which converge to the true values after $50$ days. These figures show the optimal control policies result in the system synchronization and the parameter convergence for the HIV model with time-varying parameter.

 \begin{figure}
 \centering
 \vspace*{-1.1cm}
    \includegraphics[width=13.5cm]{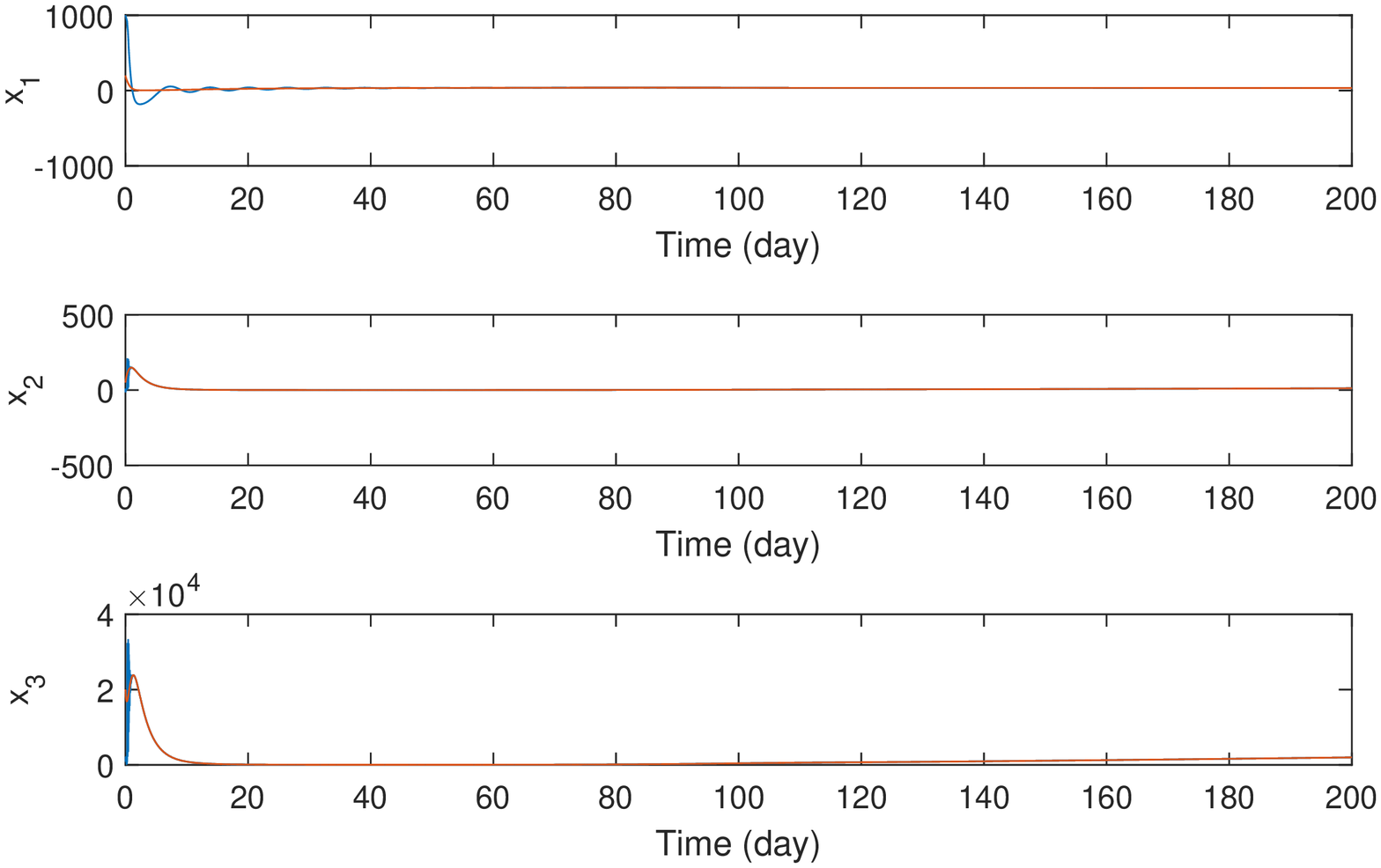}
    \vspace*{-0.5cm}
    \caption{\label{fig_sys}  (color online) Trajectories of three states of HIV model with time-varying parameter in Eq. \eqref{eq:hiv_system} where the red lines represent the dynamical evolution of drive system and the blue lines represent the evolution of response system with unknown parameters.}
 \end{figure}

\begin{figure}
 \centering
 \vspace*{-1.1cm}
    \includegraphics[width=13.5cm]{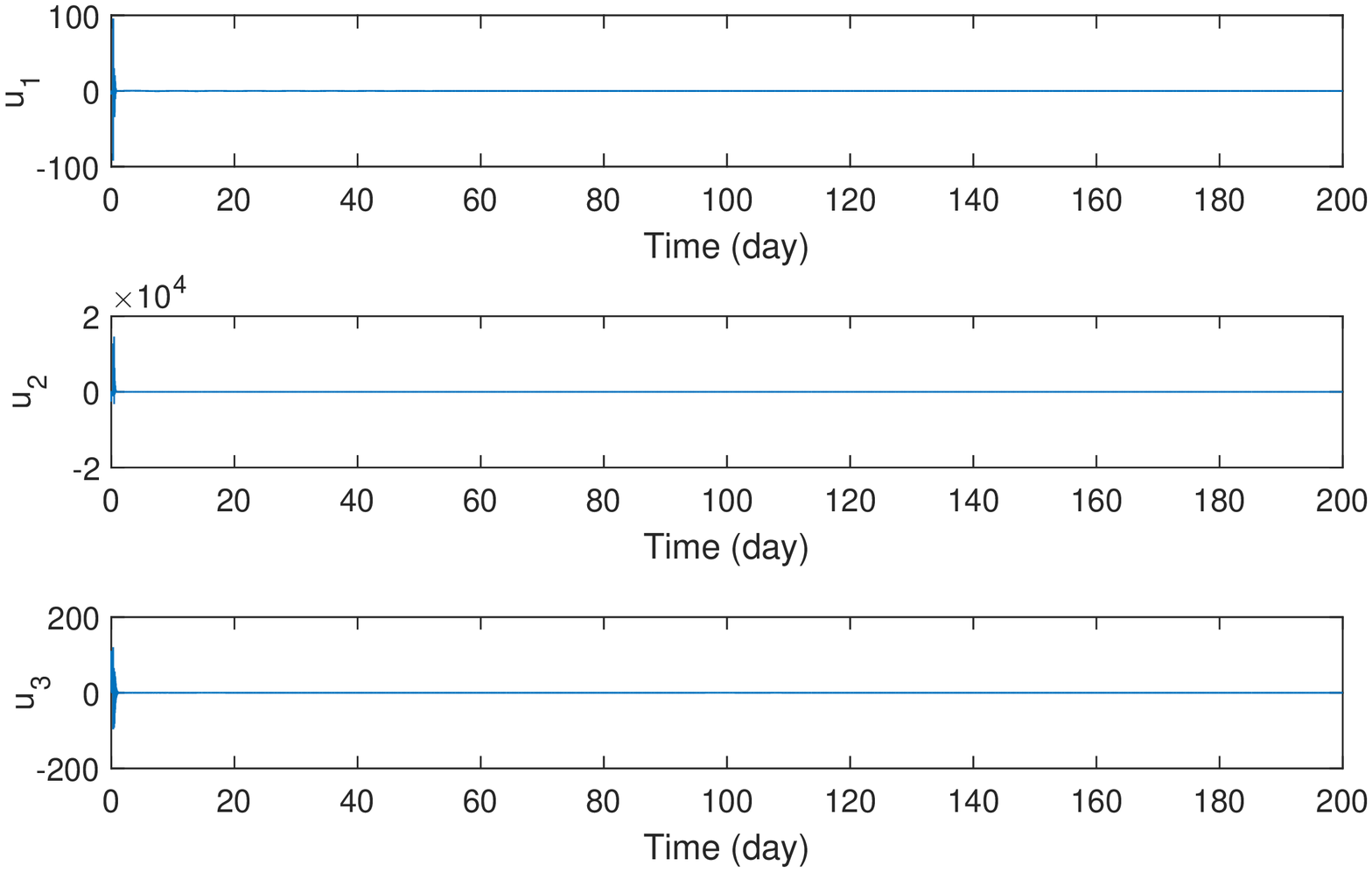}
    \vspace*{-0.5cm}
    \caption{\label{fig_u}  (color online) Trajectories of optimal control strategies generated by adaptive nonlinear receding horizon control method in Algorithm-1 for the HIV model with time-varying parameter.}
 \end{figure}
 
 \begin{figure}
 \centering
 \vspace*{-1.1cm}
    \includegraphics[width=13.5cm]{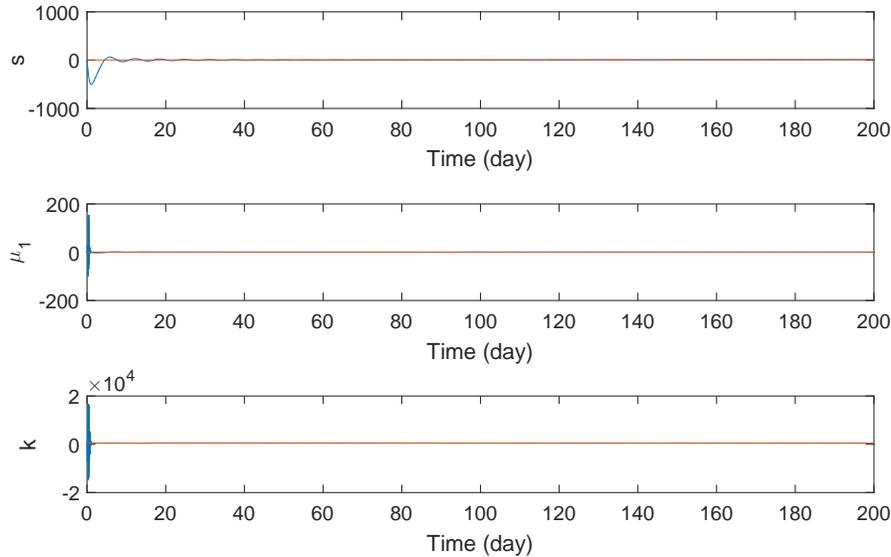}
    \vspace*{-0.5cm}
    \caption{\label{fig_par}  (color online) Trajectories of estimated dynamic parameters of HIV model with time-varying parameter where the red lines denote their true values and the blue lines denote the estimated results.}
 \end{figure}

\section{Conclusions}\label{sec:conclusions}
In this work, we proposed a real-time adaptive nonlinear receding horizon based control algorithm to estimate unknown parameters of a given nonlinear HIV system model. Under the assumption that all states are measurable, some unknown fixed and/or time-varying parameters of HIV model are shown to be identifiable, with remarkable accuracy. Compared with the existing adaptive observers for estimating parameters of HIV models (which usually need at least four or five measurements of cell count and virus load) the proposed methodology provides a \emph{one-time, single-shot, on-line algorithm} to identify the uncertain parameters, as well as generating the optimal control policy without resorting to any iterative approximation and linearlization techniques. This constitutes one of the major contributions of the study. This estimation routine can be implemented in a short period of time which enables to determine the parameters at the early infection stage, which constitutes the second major contribution of the paper. In on-going/future works, it is desired to improve the proposed method to estimate \emph{all} constant or time-varying dynamic parameters of a more complex HIV model.

\bibliography{hiv_paper}

\end{document}